\documentclass[a4paper,english,oneside]{article}
\usepackage[T1]{fontenc}
\usepackage[utf8]{inputenc} 

\usepackage{geometry}
\usepackage{amsmath}
\usepackage{amsfonts}
\usepackage{amssymb}
\usepackage{amsthm}
\usepackage[english]{babel}
\usepackage{color}
\usepackage{lmodern}
\usepackage{pstricks-add}
\usepackage{hyperref}
\hypersetup{hidelinks}

\usepackage{mathscinet}
\usepackage{enumitem}
\usepackage[cal=boondoxo,bb=ams]{mathalfa}
\usepackage{microtype}

\usepackage{mathtools}
\usepackage{esint}

\theoremstyle{plain}
\newtheorem{theorem}{Theorem}[section]
\newtheorem{lemma}[theorem]{Lemma}

\theoremstyle{definition}

\theoremstyle{remark}
\newtheorem{remark}{Remark}

\begin{document}

\title{A blow-up result for a Nakao-type weakly coupled system with nonlinearities of derivative-type}

\author{Alessandro Palmieri$\,^{\mathrm{a}}$, Hiroyuki Takamura$\,^{\mathrm{b}}$}

\date{
\small{ $\,^\mathrm{a}$ Mathematical Institute, Tohoku University, Aoba, Sendai 980-8578, Japan} \\
\small{ $\,^\mathrm{b}$ Mathematical Institute/Research Alliance Center of Mathematical Sciences, Tohoku University, Aoba, Sendai 980-8578, Japan} \\
 \normalsize{\today} }
\maketitle

\begin{abstract}
In this paper, we consider a weakly coupled system of a wave and damped Klein-Gordon equation with nonlinearities of derivative type. We prove a blow-up result for the Cauchy problem associated with this system for nonnegative and compactly supported data by means of an iteration argument.
\end{abstract}

\begin{flushleft}
\textbf{Keywords} wave equation, damped Klein-Gordon equation, iteration argument, slicing procedure, unbounded exponential multiplier, critical case
\end{flushleft}

\begin{flushleft}
\textbf{AMS Classification (2020)} 35B44, 35C15, 35L05, 35L56, 35L76
\end{flushleft}

\section{Introduction}

Let us consider a Nakao-type weakly coupled system with nonlinearities of derivative-type, namely, 
\begin{align}\label{Nakao der typ}
\begin{cases}
\partial_t^2 u-\Delta u+b\partial_t u+m^2 u=|\partial_t v|^p, & x\in\mathbb{R}^n, \, t\in (0,T), \\
\partial_t^2 v-\Delta v=|\partial_t u|^q, & x\in\mathbb{R}^n, \, t\in (0,T), \\
(u, \partial_t u)(0,x)=\varepsilon (u_0, u_1)(x), & x\in\mathbb{R}^n, \\
(v, \partial_t v)(0,x)=\varepsilon (v_0, v_1)(x), & x\in\mathbb{R}^n, 
\end{cases}
\end{align} where $p,q>1$, $\varepsilon$ is a positive parameter describing the size of the Cauchy data, and $b>0$, $m^2\geqslant 0$ are real constants.

Over the last years, systems of diffusion and wave equations with coupled nonlinear terms have been studied in the literature (see \cite{Na16,Na18,Wak17,ChenRei21,Chen21MMAS}). By diffusion equations here we mean, in a broad sense, not only parabolic equations but also hyperbolic equations which present diffusion phenomena towards certain parabolic models. This kind of nonlinear coupled systems have been named \emph{Nakao's problems} in the case of a weakly coupled Cauchy system of wave and damped wave equations in  \cite{Wak17,ChenRei21,Chen21MMAS} after the author of \cite{Na16,Na18}, who first proposed and studied these systems in the case of bounded domains.

Let us summarize briefly the results for the Nakao's problems considered in the case of the whole space, i.e. for Cauchy problems. In \cite{Wak17,ChenRei21} the Nakao's problem with weakly coupled power nonlinearities, namely,
\begin{align}\label{Nakao power nonlin}
\begin{cases}
\partial_t^2 u-\Delta u+\partial_t u=|v|^p, & x\in\mathbb{R}^n, \, t\in (0,T), \\
\partial_t^2 v-\Delta v=|u|^q, & x\in\mathbb{R}^n, \, t\in (0,T), \\
(u, \partial_t u)(0,x)=\varepsilon (u_0, u_1)(x), & x\in\mathbb{R}^n, \\
(v, \partial_t v)(0,x)=\varepsilon (v_0, v_1)(x), & x\in\mathbb{R}^n, 
\end{cases}
\end{align}
has been investigated from the viewpoint of the blow-up in finite time (for suitable $p,q$ and under suitable sign assumptions on the Cauchy data). While in \cite{Wak17} the so-called \emph{test function method} is used, in \cite{ChenRei21} an \emph{iteration argument} is employed, by considering the space averages of the components of a local solution as time-dependent functionals.

On the other hand, in \cite{Chen21MMAS} the Nakao's problem with weakly coupled nonlinearities of derivative type, namely \eqref{Nakao der typ} for $(b,m^2)=(1,0)$, is studied again from the sufficiency part. In particular, the blow-up in finite time is proved for $p,q>1$ such that $$\frac{1}{pq-1}>\frac{n-1}{2}$$ provided that the Cauchy data are compactly supported, nonnegative and nontrivial. The approach used to prove this blow-up result is inspired in some sense by \cite[Section 13.2]{LZbook} and by \cite{LT19}.

In what follows we called \eqref{Nakao der typ} a \emph{Nakao-type weakly coupled system}, since we will consider a semilinear wave equation for $v$ and a semilinear damped Klein-Gordon equation for $u$ which are weakly coupled through the nonlinear terms given by powers of the time-derivatives. We shall focus only on the case of the Cauchy problem and our goal will be determining a blow-up result in finite time when the exponents of the nonlinear terms $p,q$ belong to a suitable range and under suitable sign assumptions for the Cauchy data. 

Our approach is based on the blow-up technique introduced by Zhou in \cite{Zhou01} for the treatment of the semilinear wave equation with a nonlinearity of derivative type in all space dimensions combined with an iteration argument for determining a sequence of lower bound estimates for a suitable time-dependent functional related to a local in time solution to \eqref{Nakao der typ}. The above cite technique of Zhou consists in reducing the problem to the one-dimensional case by integrating with respect to the last $(n-1)$ space-variables and, then, in proving the blow-up on a suitable characteristic line. More specifically, when dealing with the wave equation in one space dimension, d'Alembert's formula is used to describe explicitly the solution. Consequently, before proving the main blow-up result of this paper, we are going to recall an integral representation formula for the linear equation associated with the equation for $u$ in \eqref{Nakao der typ} (which is a damped Klein-Gordon equation) in one space dimension. Moreover, since the kernel function appearing in this integral formula contain an exponential factor, we will need to adapt the treatment of an unbounded exponential multiplier in the iteration frame from \cite{ChenPal20,ChenPal20der} to our problem by applying a \emph{slicing procedure} while shrinking the domain of integration in the iteration frame. We anticipate that the other factor appearing in the integral kernel will be the composition of the modified Bessel function of the first kind of order 0 with another function related to the forward light-cone. In the derivation of the iteration frame, we will take advantage of the fact that this special function (denoted $\mathrm{I}_0$) is bounded from below by a positive function. On the contrary, we may not use the asymptotic behavior of $\mathrm{I}_0$ for large arguments due to the contemporary presence of the aforementioned exponential factor. For a rigorous explanation we address the reader to Remark \ref{Remark asymptotic}.

The range of $p,q$ for which our blow-up result is valid is exactly the same one as in \cite{Chen21MMAS} for the special case $(b,m^2)=(1,0)$ that we recalled above, although the methods employed in our proof and in the proof of the corresponding result in \cite{Chen21MMAS} are quite different. Moreover, we will extend the blow-up result even to the limit case $$\frac{1}{pq-1}=\frac{n-1}{2}.$$ 

Finally, we point out that the blow-up result in the present work is valid only under the further assumption 
\begin{align}\label{not dominant mass case}
b^2\geqslant 4 m^2.
\end{align} We refer to Remark \ref{Remark KG} for a technical explanation on the unsuitableness of our method for $b^2<4m^2$.
 We may interpret the condition \eqref{not dominant mass case} by saying that we consider the case in which the equation for $u$ in \eqref{Nakao der typ} has a \emph{mass term} $m^2 u$ that is dominated (or balanced, when the equality holds) by the \emph{damping term} $b\partial_t u$. Therefore, this equation has some properties which resemble the ones for the damped wave equation rather than the ones for the Klein-Gordon equation.
 Let us explain the previous heuristic considerations more rigorously. If we consider the linear damped Klein-Gordon equation
\begin{align*}
 \partial_t^2 \phi-\Delta \phi+b\partial_t \phi+m^2 \phi=0,
\end{align*} then, carrying out the transformation $\phi(t,x)=\mathrm{e}^{\gamma t}\psi(t,x)$, where $\gamma$ is a real constant, it results that $\psi$ solves
\begin{align*}
 \partial_t^2 \psi-\Delta \psi+(2\gamma+b)\partial_t \psi+(\gamma^2+b\gamma+m^2) \psi=0.
\end{align*} For $b^2>4m^2$ we can choose $\gamma\doteq \frac{1}{2}(-b+\sqrt{b^2-4m^2})$ so that $\psi$ solves the damped wave equation 
\begin{align*}
 \partial_t^2 \psi-\Delta \psi+(b^2-4m^2)^{\frac{1}{2}} \partial_t \psi=0.
\end{align*} For this reason, we call the case $b^2>4m^2$ the case with \emph{dominant damping}. On the contrary, for $b^2<4m^2$, setting $\gamma\doteq -\frac{b}{2}$, we get that $\psi$ solves the Klein-Gordon equation (with positive mass)
\begin{align*}
 \partial_t^2 \psi-\Delta \psi+\left(m^2-\frac{b^2}{4}\right) \psi=0.
\end{align*} Hence, we call  $b^2<4m^2$ the case with \emph{dominant mass}. In the limit case $b^2=4m^2$, we find that $\psi$ solves the free wave equation, therefore, we call it the \emph{balanced case}. We stress that this nomenclature is borrowed from the introduction of \cite{ER18}.

The paper is organized as follows: in Section \ref{Section main thm} we state the main blow-up result for \eqref{Nakao der typ}; in Section \ref{Section repres formula} we recall the integral representation formula for the linear Cauchy problem associated with the damped Klein-Gordon equation when $n=1$; finally, in Section \ref{Section proof main thm} we derive the iteration frame 
and we apply the slicing procedure to perform the iteration procedure.

\section{Main result} \label{Section main thm}

\begin{theorem} \label{Main Thm} Let $n\geqslant 1$ and let $b>0,m^2\geqslant 0$ be real constants satisfying \eqref{not dominant mass case}. We assume that $u_0,v_0\in \mathcal{C}^2_0(\mathbb{R}^n)$,  $u_1,v_1\in \mathcal{C}^1_0(\mathbb{R}^n)$ are nonnegative and compactly supported functions with supports contained in $B_R$ for some $R>0$, and that $v_1$ is nontrivial. Let us consider exponents for the nonlinear terms $p,q>1$ satisfying 
\begin{align}\label{def theta}
\theta(n,p,q)\doteq \frac{1}{pq-1} -\frac{n-1}{2}\geqslant 0.
\end{align} Then, there exists $\varepsilon_0=\varepsilon_0(n,p,q,b,R,v_1)$ such that for any $\varepsilon\in (0,\varepsilon_0]$ if $(u,v)\in \big(\mathcal{C}([0,T)\times \mathbb{R}^n)\big)^2$ is a local in time solution to \eqref{Nakao der typ} such that
\begin{align}\label{support condition}
\mathrm{supp} \, u(t,\cdot), \mathrm{supp} \, v(t,\cdot) \subset B_{R+t} \quad \mbox{for any} \ t\in [0,T),
\end{align} where $T=T(\varepsilon)$ denotes the lifespan of $(u,v)$, then, $(u,v)$ blows up in finite time.

Furthermore, the following upper bound estimate for the lifespan holds
\begin{align}\label{upper bound lifespan thm}
T(\varepsilon)\leqslant \begin{cases} C \varepsilon^{-\theta(n,p,q)^{-1}} & \mbox{if} \ \ \theta(n,p,q)>0, \\
\exp\left(C\varepsilon^{-(pq-1)}\right) & \mbox{if} \ \ \theta(n,p,q)=0,
\end{cases}
\end{align}
 where the positive constant $C$ is independent of $\varepsilon$.
\end{theorem}

\section{Integral representation formula in one space dimension} \label{Section repres formula}

In the proof of Theorem \ref{Main Thm} we are going to use the approach from \cite{Zhou01} to proving the blow-up on a certain characteristic line, as described in the introduction. 

Since the second order partial operator acting on $u$ in \eqref{Nakao der typ} is a damped wave operator with a mass term we need first to get a representation formula for the corresponding linear Cauchy problem in the one-dimensional case, namely,
\begin{align}\label{CP lin b,m 1D}
\begin{cases}
\partial_t^2 \phi -\partial_x^2 \phi+b\partial_t \phi+m^2\phi=F(t,x), & x\in\mathbb{R}, \, t>0, \\
\phi(0,x)=f(x), & x\in\mathbb{R}, \\
\partial_t \phi(0,x)=g(x), & x\in\mathbb{R}. \\
\end{cases}
\end{align}

The integral representation formula for the solution to \eqref{CP lin b,m 1D}, under suitable regularity assumptions on the data $f,g,F$, is already known in the literature. However, the proof of this representation formula in the form that we will employ is scattered through different references. For the ease of readability we shall provide an elementary proof of it. 

In what follows, we collect and adapt the results from \cite[Chapter III Section 3.5 and Chapter VI Section 12.6]{CHvol2} and \cite[Section 1.1]{Yag10}.

\begin{lemma}\label{Lemma RF} Let $b>0$ and $m^2\geqslant 0$. For any $h\in\mathcal{C}^1(\mathbb{R})$ and any $t\geqslant 0, x\in\mathbb{R}$ we define the solution operator
\begin{align} \label{sol op}
\mathrm{S}\big(t;b,m^2\big) h(x)\doteq 
\begin{cases} \displaystyle{\frac{1}{2}  \, \mathrm{e}^{-\frac{b}{2} t}  \int_{x-t}^{x+t} \mathrm{I}_0\left(\mu\sqrt{t^2-|x-y|^2} \,\right)h(y) \, \mathrm{d}y} & \mbox{for} \ \, 4m^2< b^2, \\
\displaystyle{\frac{1}{2}  \, \mathrm{e}^{-\frac{b}{2}t}  \int_{x-t}^{x+t} h(y) \, \mathrm{d}y} & \mbox{for} \ \, 4m^2= b^2, \phantom{\Bigg(}\\
\displaystyle{\frac{1}{2}  \, \mathrm{e}^{-\frac{b}{2}t}  \int_{x-t}^{x+t} \mathrm{J}_0\left(\mu\sqrt{t^2-|x-y|^2} \,\right)h(y) \, \mathrm{d}y} & \mbox{for} \ \, 4m^2> b^2, 
\end{cases}
\end{align}  where $$\mu\doteq \sqrt{\left|\frac{b^2}{4}-m^2\right|}$$ and $\mathrm{I}_0,\mathrm{J}_0$ denote the modified Bessel function and the Bessel function of the first kind of order $0$, respectively, (cf. \cite[Sections 10.2 and 10.25]{OLBC10}).

Let us consider $f\in\mathcal{C}^2(\mathbb{R}),g\in\mathcal{C}^1(\mathbb{R})$ and $F\in \mathcal{C}^1([0,\infty)\times \mathbb{R})$. Then, the solution to the linear Cauchy problem \eqref{CP lin b,m 1D} is given by
\begin{align}\label{representation formula}
\phi(t,x)= \mathrm{S}\big(t;b,m^2\big)(g+bf)(x)+\frac{\partial}{\partial t} \mathrm{S}\big(t;b,m^2\big)f(x)+\int_0^t \mathrm{S}\big(t-\tau;b,m^2\big)(F(\tau,\cdot))(x)\, \mathrm{d}\tau.
\end{align}
\end{lemma}

\begin{remark}
In the special case $(b,m^2)=(1,0)$ the representation formula \eqref{representation formula} coincides with the one for the classical linear damped wave equation (see \cite[Equation (43), page 695]{CHvol2} or \cite[Proposition 2.1]{SW17}).
\end{remark}

\begin{proof}
 In the balanced case $b^2=4m^2$ the function $\psi(t,x)=\mathrm{e}^{\frac{b}{2}t}\phi(t,x)$ solves the Cauchy problem
\begin{align*}
\begin{cases}
\partial_t^2 \psi -\partial_x^2 \psi=\mathrm{e}^{\frac{b}{2}t}F(t,x), & x\in\mathbb{R}, \, t>0, \\
\phi(0,x)=f(x), & x\in\mathbb{R}, \\
\partial_t \phi(0,x)=g(x)+\frac{b}{2}f(x), & x\in\mathbb{R}. \\
\end{cases}
\end{align*} Combining d'Alembert's formula with Duhamel's principle and the inverse transformation $\phi(t,x)=\mathrm{e}^{-\frac{b}{2}t}\psi(t,x)$, we get immediately \eqref{representation formula}.

When $b^2\neq 4m^2$ we begin by proving that $\mathrm{S}\big(t;b,m^2\big)(g)(x)$ solves the Cauchy problem \eqref{CP lin b,m 1D} for $f=0$ and $F=0$. We carry on the computation only in the dominant damping case $b^2>4m^2$, since in the dominant mass case $b^2<4m^2$ the procedure is completely analogous. Let us check the Cauchy conditions first. Clearly $\mathrm{S}\big(0;b,m^2\big)(g)(x)=0$. On the other hand, using $\mathrm{I}_0(0)=1$, we have
\begin{align}
\frac{\partial}{\partial t}\mathrm{S}\big(t;b,m^2\big)(g)(x) & = \frac{1}{2} \, \mathrm{e}^{-\frac{b}{2}t}(g(x+t)+g(x-t))-\frac{b}{4}  \, \mathrm{e}^{-\frac{b}{2} t}  \int_{x-t}^{x+t} \mathrm{I}_0\left(\mu\sqrt{t^2-|x-y|^2} \,\right)g(y) \, \mathrm{d}y \notag\\ & \qquad +\frac{\mu}{2} \, t  \, \mathrm{e}^{-\frac{b}{2} t}  \int_{x-t}^{x+t} \frac{\mathrm{I}_0'\left(\mu\sqrt{t^2-|x-y|^2}  \right)}{\sqrt{t^2-|x-y|^2} } g(y) \, \mathrm{d}y. \label{time derivative S(t)g}
\end{align} Consequently, $\partial_t\mathrm{S}\big(t;b,m^2\big)(g)(x)\big|_{t=0}=g(x)$.

We prove now that $\mathrm{S}\big(t;b,m^2\big)(g)(x)$ solves the homogeneous differential equation. A further differentiation of \eqref{time derivative S(t)g} with respect to $t$ provides
\begin{align} 
& \frac{\partial^2}{\partial t^2}\mathrm{S}\big(t;b,m^2\big)(g)(x) \notag \\ 
& \quad  = \left(-\frac{b}{2} +\frac{\mu^2 t}{4} \right) \mathrm{e}^{-\frac{b}{2}t}(g(x+t)+g(x-t))+\frac{1}{2} \, \mathrm{e}^{-\frac{b}{2}t}(g'(x+t)-g'(x-t)) \notag \\
 & \qquad +\frac{b^2}{8}  \, \mathrm{e}^{-\frac{b}{2} t}  \int_{x-t}^{x+t} \mathrm{I}_0\left(\mu\sqrt{t^2-|x-y|^2} \,\right) g(y) \, \mathrm{d}y \notag\\
  & \qquad +\frac{1}{2}  \, \mathrm{e}^{-\frac{b}{2} t}  \int_{x-t}^{x+t} \mathrm{I}_0'\left(\mu\sqrt{t^2-|x-y|^2}  \right)\left(\frac{\mu(1-bt)}{\sqrt{t^2-|x-y|^2} } -\frac{\mu t^2}{(t^2-|x-y|^2)^{3/2}}\right) g(y) \, \mathrm{d}y \notag \\ 
 & \qquad  +\frac{1}{2}  \, \mathrm{e}^{-\frac{b}{2} t}  \int_{x-t}^{x+t} \mathrm{I}_0''\left(\mu\sqrt{t^2-|x-y|^2}  \right)\frac{\mu^2t^2}{t^2-|x-y|^2 } \, g(y) \, \mathrm{d}y.\label{second time derivative S(t)g}
\end{align} We point out that, differentiating the second integral in \eqref{time derivative S(t)g}, we applied the relation 
\begin{align}\label{I1(z)/z for z=0}
\frac{\mathrm{I}_0'(z)}{z}\bigg|_{z=0}=\frac{1}{2}
\end{align} that follows from the relation $\mathrm{I}_0'=\mathrm{I}_1$ and from the Maclaurin series expansion for the function $z^{-1}\mathrm{I}_1(z)$ (cf. \cite[Equations  (10.29.3) and (10.25.2)]{OLBC10}). Using again \eqref{I1(z)/z for z=0}, we find that the second order derivative with respect to $x$ of $\mathrm{S}\big(t;b,m^2\big)(g)(x)$ is given by
\begin{align}
& \frac{\partial^2}{\partial x^2}\mathrm{S}\big(t;b,m^2\big)(g)(x) \notag \\ 
& \quad  = \frac{\mu^2 t}{4} \mathrm{e}^{-\frac{b}{2}t}(g(x+t)+g(x-t))+\frac{1}{2} \, \mathrm{e}^{-\frac{b}{2}t}(g'(x+t)-g'(x-t)) \notag \\
  & \qquad +\frac{1}{2}  \, \mathrm{e}^{-\frac{b}{2} t}  \int_{x-t}^{x+t} \mathrm{I}_0'\left(\mu\sqrt{t^2-|x-y|^2}  \right)\left(-\frac{\mu}{\sqrt{t^2-|x-y|^2} } -\frac{\mu (x-y)^2}{(t^2-|x-y|^2)^{3/2}}\right) g(y) \, \mathrm{d}y \notag \\ 
 & \qquad  +\frac{1}{2}  \, \mathrm{e}^{-\frac{b}{2} t}  \int_{x-t}^{x+t} \mathrm{I}_0''\left(\mu\sqrt{t^2-|x-y|^2}  \right)\frac{\mu^2(x-y)^2}{t^2-|x-y|^2 } \, g(y) \, \mathrm{d}y.
\label{second space derivative S(t)g}
\end{align} 
Combining \eqref{time derivative S(t)g}, \eqref{second time derivative S(t)g} and \eqref{second space derivative S(t)g}, we get
\begin{align}
& \left(\frac{\partial^2}{\partial t^2}-\frac{\partial^2}{\partial x^2}+b\frac{\partial}{\partial t}+m^2 I\right)  \mathrm{S}\big(t;b,m^2\big)(g)(x) \notag \\
& \quad =  \frac{1}{2}  \, \mathrm{e}^{-\frac{b}{2} t}  \int_{x-t}^{x+t}  \mu^2 \mathrm{I}_0''\left(\mu\sqrt{t^2-|x-y|^2 }\right)  g(y) \, \mathrm{d}y + \frac{1}{2}  \, \mathrm{e}^{-\frac{b}{2} t}  \int_{x-t}^{x+t}  \frac{\mu  \mathrm{I}_0'\left(\mu\sqrt{t^2-|x-y|^2} \right)}{\sqrt{t^2-|x-y|^2}}g(y) \, \mathrm{d}y \notag \\ & \qquad +\frac{1}{2}  \, \mathrm{e}^{-\frac{b}{2} t}  \int_{x-t}^{x+t} \left(m^2-\frac{b^2}{4}\right) \mathrm{I}_0\left(\mu\sqrt{t^2-|x-y|^2} \right)  g(y) \, \mathrm{d}y \notag \\
& \quad =  \frac{\mu^2}{2}  \, \mathrm{e}^{-\frac{b}{2} t}  \int_{x-t}^{x+t} \left(\mathrm{I}_0''(z)+\frac{\mathrm{I}'_0(z)}{z}-\mathrm{I}_0(z)\right)\bigg|_{z=\mu\sqrt{t^2-(x-y)^2}} \ g(y) \, \mathrm{d}y =0, \label{eq S(t)g}
\end{align} where in the last step we used the fact that $\mathrm{I}_0$ is a solution of the ODE (see \cite[Equation (10.25.1)]{OLBC10})
 $$z^2\mathrm{I}_0''(z)+z \mathrm{I}'_0(z)-z^2\mathrm{I}_0(z)=0.$$
We emphasize that in the dominant mass case we can repeat the same steps as before. However, since $\mu^2=m^2-\frac{b^2}{4}$ in this case, we use the fact that $\mathrm{J}_0$ is a solution of the ODE (see \cite[Equation (10.2.1)]{OLBC10})
 $$z^2\mathrm{J}_0''(z)+z \mathrm{J}'_0(z)+z^2\mathrm{J}_0(z)=0.$$
So, we proved \eqref{representation formula} for $f=0$ and $F=0$.

 Now we focus on the case $g=0$ and $F=0$. We claim that $$\widetilde{\phi}(t,x)\doteq \frac{\partial}{\partial t}\mathrm{S}\big(t;b,m^2\big)(f)(x)+b\mathrm{S}\big(t;b,m^2\big)(f)(x)$$ is the solution of \eqref{CP lin b,m 1D} with vanishing second data and source term. 

Clearly, $\widetilde{\phi}$ solves the homogeneous differential equation as the differential operators $(\partial_t+bI)$ and $(\partial_t^2-\partial_x^2+b\partial_t+m^2 I)$ commute. We check now the Cauchy conditions. Using the initial conditions derived in the previous case, we see immediately that $\widetilde{\phi}(0,x)=f(x)$. On the other hand,
\begin{align*}
\frac{\partial}{\partial t}\widetilde{\phi}(t,x)= \left(\frac{\partial^2}{\partial t^2}+b\frac{\partial}{\partial t}\right)\mathrm{S}\big(t;b,m^2\big)(f)(x) = \left(\frac{\partial^2}{\partial x^2}-m^2 I\right)\mathrm{S}\big(t;b,m^2\big)(f)(x). 
\end{align*} Therefore, combining \eqref{sol op} and \eqref{second space derivative S(t)g} with the previous relation it follows that $\partial_t \widetilde{\phi}(0,x)=0$.

It remains to consider the inhomogeneous Cauchy problem \eqref{CP lin b,m 1D} with both vanishing initial data $f=g=0$. By using Duhamel's principle together with the solution operator defined in \eqref{sol op}, since the model under consideration is invariant by time translations, we get that the solution for this case is given by $$\int_0^t \mathrm{S}\big(t-\tau;b,m^2\big)(F(\tau,\cdot))(x)\, \mathrm{d}\tau.$$ 
Due to the linearity of \eqref{CP lin b,m 1D}, combining the results from the previous subcases, we conclude the validity of \eqref{representation formula}.
\end{proof}


\begin{remark} By using \eqref{sol op} and \eqref{time derivative S(t)g}, we can rewrite \eqref{representation formula} more explicitly as follows:
\begin{align}
\phi(t,x) &= \frac{1}{2} \, \mathrm{e}^{-\frac{b}{2}t}(f(x+t)+f(x-t))+\frac{1}{2}  \, \mathrm{e}^{-\frac{b}{2} t}  \int_{x-t}^{x+t} \mathrm{I}_0\left(\mu\sqrt{t^2-|x-y|^2} \,\right) \left(g(y)+\frac b 2 f(y)\right) \, \mathrm{d}y \notag\\ & \qquad +\frac{\mu}{2} \, t  \, \mathrm{e}^{-\frac{b}{2} t}  \int_{x-t}^{x+t} \frac{\mathrm{I}_1\left(\mu\sqrt{t^2-|x-y|^2}  \right)}{\sqrt{t^2-|x-y|^2} } f(y) \, \mathrm{d}y \notag\\ 
& \qquad + \frac{1}{2} \int_0^t \mathrm{e}^{-\frac{b}{2}(t-\tau)}\int_{x-t+\tau}^{x+t-\tau}  \mathrm{I}_0\left(\mu\sqrt{(t-\tau)^2-|x-y|^2} \,\right) F(\tau,y) \, \mathrm{d}y \, \mathrm{d}\tau
\label{representation formula dom damp}
\end{align}
for $b^2>4m^2$, and 
\begin{align}
\phi(t,x) &= \frac{1}{2} \, \mathrm{e}^{-\frac{b}{2}t}(f(x+t)+f(x-t))+\frac{1}{2}  \, \mathrm{e}^{-\frac{b}{2} t}  \int_{x-t}^{x+t} \left(g(y)+\frac b 2 f(y)\right) \, \mathrm{d}y \notag\\ 
& \qquad + \frac{1}{2} \int_0^t \mathrm{e}^{-\frac{b}{2}(t-\tau)}\int_{x-t+\tau}^{x+t-\tau}  F(\tau,y) \, \mathrm{d}y \, \mathrm{d}\tau
\label{representation formula balanced case}
\end{align}
for $b^2=4m^2$.

Finally, for  $b^2<4m^2$ the representation formula is analogous the the one in \eqref{representation formula dom damp}, but instead of the modified Bessel functions $\mathrm{I}_0,\mathrm{I}_1$ we have the Bessel functions $\mathrm{J}_0,-\mathrm{J}_1$, respectively. In particular, we use the relation $\mathrm{J}_0'=-\mathrm{J}_1$, see \cite[Equation (10.6.2)]{OLBC10}.
\end{remark}

\begin{remark} \label{Remark KG} In the statement of Theorem \ref{Main Thm} we consider only $b,m^2$ such that $b^2\geqslant 4m^2$. This assumption is due to the fact in the dominant mass case $b^2<4m^2$ the kernel functions in the representation formula \eqref{representation formula} are no longer nonnegative functions. Indeed, in the iteration argument that we will use to prove the blow-up result it is crucial the fact that we will be working with a nonnegative functional. For $b^2<4m^2$ the partial differential operator acting on $u$ in \eqref{Nakao der typ} is in this sense very close to the Klein-Gordon operator (i.e. for $b=0$) and the damped oscillations of the Bessel functions of the first kind do not allow to carry on with the iteration procedure.
\end{remark}

\section{Proof of Theorem \ref{Main Thm}} \label{Section proof main thm}

The proof of Theorem \ref{Main Thm} is based on the approach introduced by Zhou in \cite{Zhou01}, where a blow-up result for the semilinear wave equation with nonlinearity of derivative-type is proved for all space dimensions. Recently, this approach have been applied  
to study semilinear models with time-dependent coefficients (cf. \cite{PalTu21,LP21,HHP21}). 

In \cite{Zhou01} d'Alembert's formula is used to prove the blow-up result for the semilinear wave equation with nonlinearity of derivative type. In our case, since we work with the weakly coupled system \eqref{Nakao der typ} together with d'Alembert's formula (coming from the equation for $v$) we shall also employ the representation formulas \eqref{representation formula dom damp} and \eqref{representation formula balanced case} from Section \ref{Section repres formula}. Notice that \eqref{representation formula dom damp} coincides exactly with \eqref{representation formula balanced case} for $\mu=0$. Hence, in what follows we work always with \eqref{representation formula dom damp} for both cases.

Let us introduce the following notation: we will write any $x\in\mathbb{R}$ as $x=(z,w)$ with $z\in\mathbb{R}$ and $w\in\mathbb{R}^{n-1}$. Thanks to this notation we might introduce the following functions
\begin{align*}
&\mathcal{U}(t,z)\doteq \int_{\mathbb{R}^{n-1}}u(t,z,w) \, \mathrm{d}w, \quad \mathcal{V}(t,z)\doteq \int_{\mathbb{R}^{n-1}}v(t,z,w) \, \mathrm{d}w \quad \mbox{for any} \ t\in [0,T), \, z\in\mathbb{R}, \\
&\mathcal{U}_j(z)\doteq \int_{\mathbb{R}^{n-1}}u_j(z,w) \, \mathrm{d}w,  \quad \ \ \mathcal{V}_j(z)\doteq \int_{\mathbb{R}^{n-1}}v_j(z,w) \, \mathrm{d}w \qquad \, \mbox{for any} \  z\in\mathbb{R}, \, j=0,1. 
\end{align*}
Clearly, it makes sense to introduce these functions only for $n\geqslant 2$, while for $n=1$ we set simply $(\mathcal{U},\mathcal{V})=(u,v)$ and $(\mathcal{U}_0,\mathcal{U}_1,\mathcal{V}_0,\mathcal{V}_1)=(u_0,u_1,v_0,v_1)$.

We remark that due to the assumption $\mathrm{supp}\, u_j, \mathrm{supp} \, v_j \subset B_R$ for $j=0,1$ it follows that
\begin{align}\label{support condition aux fun CD}
\mathrm{supp}\, \mathcal{U}_j , \mathrm{supp}\, \mathcal{V}_j \subset (-R,R) \quad j=0,1.
\end{align} Analogously, from \eqref{support condition} we have
\begin{align} \label{support condition aux fun}
\mathrm{supp}\, \mathcal{U}(t,\cdot) , \mathrm{supp}\, \mathcal{V}(t,\cdot) \subset \big(-(R+t),R+t\big) \quad \mbox{for any} \ t\in[0,T). 
\end{align}
By a straightforward computation we find that  $(\mathcal{U},\mathcal{V})$ solves for $n\geqslant 2$ the following system
\begin{align*}
\begin{cases}
\partial_t^2 \mathcal{U}-\partial_z^2 \mathcal{U}+b \partial_t \mathcal{U}+ m^2 \mathcal{U}=\displaystyle{\int_{\mathbb{R}^{n-1}}|\partial_t v(t,z,w)|^p \mathrm{d}w}, &  t\in (0,T), \, z\in\mathbb{R}^n, \\
\partial_t^2 \mathcal{V}-\partial_z^2 \mathcal{V}=\displaystyle{\int_{\mathbb{R}^{n-1}}|\partial_t u(t,z,w)|^q \mathrm{d}w}, & t\in (0,T), \, z\in\mathbb{R}^n, \\
(\mathcal{U}, \partial_t \mathcal{U})(0,z)=\varepsilon (\mathcal{U}_0, \mathcal{U}_1)(z), & z\in\mathbb{R}^n, \\
(\mathcal{V}, \partial_t \mathcal{V})(0,z)=\varepsilon (\mathcal{V}_0, \mathcal{V}_1)(z), & z\in\mathbb{R}^n.
\end{cases}
\end{align*}
By using D'Alembert's formula and the representation formula for the damped wave equation with a mass term from Section \ref{Section repres formula}, we obtain the following integral representations
\begin{align*}
\mathcal{U}(t,z) &= \mathcal{U}^{\mathrm{lin}}(t,z)+\mathcal{U}^{\mathrm{nlin}}(t,z), \\
\mathcal{V}(t,z) &= \mathcal{V}^{\mathrm{lin}}(t,z)+\mathcal{V}^{\mathrm{nlin}}(t,z), 
\end{align*}
where
\begin{align*}
\mathcal{U}^{\mathrm{lin}}(t,z) &\doteq \frac{\varepsilon}{2}\,\mathrm{e}^{-\frac{b}{2}t} \big(\mathcal{U}_0(z+t)+\mathcal{U}_0(z-t)\big)+ \frac{\varepsilon}{2} \,\mathrm{e}^{-\frac{b}{2}t} \int_{z-t}^{z+t}\mathrm{I}_0\left(\mu\sqrt{t^2-|z-y|^2} \,\right)\big(\mathcal{U}_1(y)+\tfrac{b}{2} \mathcal{U}_0(y)\big) \, \mathrm{d}y \\
& \quad +\frac{\mu \, \varepsilon}{2} \, t\, \mathrm{e}^{-\frac{b}{2}t} \int_{z-t}^{z+t} \frac{\mathrm{I}_1\left(\mu\sqrt{t^2-|z-y|^2} \,\right)}{\sqrt{t^2-|z-y|^2}} \, \mathcal{U}_0(y) \, \mathrm{d}y, \\
\mathcal{U}^{\mathrm{nlin}}(t,z) &\doteq \frac{1}{2} \int_0^t \mathrm{e}^{-\frac{b}{2}(t-\tau)}\int_{z-t+\tau}^{z+t-\tau} \mathrm{I}_0\left(\mu\sqrt{(t-\tau)^2-|z-y|^2} \,\right) \int_{\mathbb{R}^{n-1}}|\partial_t v(\tau,y,w)|^p \mathrm{d}w \, \mathrm{d}y \, \mathrm{d}\tau, \\
\mathcal{V}^{\mathrm{lin}}(t,z) &\doteq \frac{\varepsilon}{2} \big(\mathcal{V}_0(z+t)+\mathcal{V}_0(z-t)\big)+ \frac{\varepsilon}{2} \int_{z-t}^{z+t}\mathcal{V}_1(y) \, \mathrm{d}y, \\
\mathcal{V}^{\mathrm{nlin}}(t,z) &\doteq \frac{1}{2} \int_0^t \int_{z-t+\tau}^{z+t-\tau} \int_{\mathbb{R}^{n-1}}|\partial_t u(\tau,y,w)|^q \mathrm{d}w \, \mathrm{d}y \, \mathrm{d}\tau. 
\end{align*}

Now that we obtained the explicit integral representation formulas for $(\mathcal{U},\mathcal{V})$, we need to determine the functional related to $(u,v)$ that blows up in finite time. We anticipate that this functional will be $\mathcal{V}$ evaluated on a certain characteristic line. In order to prove the blow-up result we will establish a sequence of lower bound estimates for this functional, that we will determine by means of a suitable iteration frame.

The next step is to determine the iteration frame. For this purpose we proceed with lower bound estimates for the functions $\mathcal{U}^{\mathrm{nlin}},\mathcal{V}^{\mathrm{nlin}}$. Hereafter we focus on the case $n\geqslant 2$, nevertheless our computations can be repeated with simple modifications in the case $n=1$.

By the support condition \eqref{support condition} we get 
\begin{align*}
\mathrm{supp} \, \partial_t u(t,\cdot), \,  \mathrm{supp} \, \partial_t  v(t,\cdot) \subset B_{R+t} \quad \mbox{for any} \ t\in [0,T),
\end{align*} that implies in turn
\begin{align}\label{support condition t der wrt t,z}
\mathrm{supp} \, \partial_t u(t,z,\cdot), \,  \mathrm{supp} \, \partial_t  v(t,z,\cdot) \subset \left\{w\in\mathbb{R}^{n-1}: |w|\leqslant\left((R+t)^2-z^2\right)^{1/2}\right\}
\end{align} for any $t\in [0,T)$ and any $z\in\mathbb{R}$ such that $|z|\leqslant R+t$.
Combining H\"older's inequality and \eqref{support condition t der wrt t,z}, we arrive at
\begin{align*}
& \int_{\mathbb{R}^{n-1}}|\partial_t v(\tau,y,w)|^p \mathrm{d}w \gtrsim \big((R+\tau)^2-y^2\big)^{-\frac{n-1}{2}(p-1)} |\partial_t \mathcal{V}(\tau,y)|^p, \\
& \int_{\mathbb{R}^{n-1}}|\partial_t u(\tau,y,w)|^q \mathrm{d}w \gtrsim \big((R+\tau)^2-y^2\big)^{-\frac{n-1}{2}(q-1)} |\partial_t \mathcal{U}(\tau,y)|^q,
\end{align*} for any $\tau\in[0,t]$ and any $y\in [z-t+\tau,z+t-\tau]$. Thus, we obtain
\begin{align*}
& \mathcal{U}^{\mathrm{nlin}}(t,z) \gtrsim  \int_0^t \mathrm{e}^{-\frac{b}{2}(t-\tau)}\int_{z-t+\tau}^{z+t-\tau} \mathrm{I}_0\left(\mu\sqrt{(t-\tau)^2-|z-y|^2} \,\right) \big((R+\tau)^2-y^2\big)^{-\frac{n-1}{2}(p-1)} |\partial_t \mathcal{V}(\tau,y)|^p \,\mathrm{d}y \,\mathrm{d}\tau,  \\
&\mathcal{V}^{\mathrm{nlin}}(t,z) \gtrsim  \int_0^t \int_{z-t+\tau}^{z+t-\tau} \big((R+\tau)^2-y^2\big)^{-\frac{n-1}{2}(q-1)} |\partial_t \mathcal{U}(\tau,y)|^q\, \mathrm{d}y \, \mathrm{d}\tau.
\end{align*} Applying Fubini's theorem, we have
\begin{align*}
& \mathcal{U}^{\mathrm{nlin}}(t,z) \gtrsim  \int_{z-t}^{z+t}  \int_0^{t-|z-y|}  \mathrm{e}^{-\frac{b}{2}(t-\tau)} \mathrm{I}_0\left(\mu\sqrt{(t-\tau)^2-|z-y|^2} \,\right) \big((R+\tau)^2-y^2\big)^{-\frac{n-1}{2}(p-1)} |\partial_t \mathcal{V}(\tau,y)|^p \, \mathrm{d}\tau \, \mathrm{d}y,  \\
&\mathcal{V}^{\mathrm{nlin}}(t,z) \gtrsim \int_{z-t}^{z+t}  \int_0^{t-|z-y|} \big((R+\tau)^2-y^2\big)^{-\frac{n-1}{2}(q-1)} |\partial_t \mathcal{U}(\tau,y)|^q  \, \mathrm{d}\tau \, \mathrm{d}y.
\end{align*}
From here on we will work on the characteristic line $t-z=R$ for $z\geqslant R$. Also, shrinking the domain of integration in the previous estimate for $\mathcal{U}^{\mathrm{nlin}}$, we find
\begin{align*}
\mathcal{U}^{\mathrm{nlin}}(R+z,z) & \gtrsim \int_{R}^{z}  \int_{y-R}^{y+R}  \mathrm{e}^{-\frac{b}{2}(t-\tau)} \mathrm{I}_0\left(\mu \sqrt{(t-\tau)^2-|z-y|^2} \,\right) \big((R+\tau)^2-y^2\big)^{-\frac{n-1}{2}(p-1)} |\partial_t \mathcal{V}(\tau,y)|^p \, \mathrm{d}\tau \, \mathrm{d}y \\
& \gtrsim \int_{R}^{z} \mathrm{e}^{-\frac{b}{2}(z-y)} (R+y)^{-\frac{n-1}{2}(p-1)}\int_{y-R}^{y+R}   \mathrm{I}_0\left(\mu\sqrt{(t-\tau)^2-|z-y|^2} \,\right)  |\partial_t \mathcal{V}(\tau,y)|^p \, \mathrm{d}\tau \, \mathrm{d}y \\
& \gtrsim \int_{R}^{z} \mathrm{e}^{-\frac{b}{2}(z-y)} (R+y)^{-\frac{n-1}{2}(p-1)}\int_{y-R}^{y+R}  |\partial_t \mathcal{V}(\tau,y)|^p \, \mathrm{d}\tau \, \mathrm{d}y,
\end{align*} where in the last step we used the inequality $\mathrm{I}_0(s)\geqslant 1$ for any $s\geqslant 0$ (due to $\mathrm{I}_0'(s)=\mathrm{I}_1(s)\geqslant 0$ for any $s\geqslant 0$ and $\mathrm{I}_0(0)=1$).

\begin{remark} \label{Remark asymptotic} As we pointed out in the introduction, we may not use the asymptotic estimate $$\mathrm{I}_0(s)\sim \frac{1}{\sqrt{2\pi s}} \mathrm{e}^s \quad \mbox{for} \  s\to \infty$$ while deriving the previous inequality. Indeed, on the domain of integration (namely, for $y\in [R,z]$ and $\tau\in [y-R,y+R]$) the argument of the modified Bessel function of the first kind of order 0 satisfies $$\mu \sqrt{(t-\tau)^2-|z-y|^2}\approx \mu \sqrt{z-y},$$ so it can be large only for $y$ away from a neighborhood of $z$. However, if we shrink further the domain of integration by removing a neighborhood of $z$, then, we are not able to compensate the exponentially decaying term $\mathrm{e}^{-\frac{b}{2}z}$ through the factor $\mathrm{e}^{\frac{b}{2}y}$ in the integral. This explains why earlier we had to use the lower bound estimate $\mathrm{I}_0\geqslant 1$ rather than the asymptotic estimate for $\mathrm{I}_0$.
\end{remark}

Then, by Jensen's inequality and the fundamental theorem of calculus, we get
\begin{align}
\mathcal{U}^{\mathrm{nlin}}(R+z,z) &  \gtrsim \int_{R}^{z} \mathrm{e}^{-\frac{b}{2}(z-y)} (R+y)^{-\frac{n-1}{2}(p-1)}\bigg| \int_{y-R}^{y+R}  \partial_t \mathcal{V}(\tau,y) \, \mathrm{d}\tau  \, \bigg|^p \mathrm{d}y \notag \\
&  \gtrsim \int_{R}^{z} \mathrm{e}^{-\frac{b}{2}(z-y)} (R+y)^{-\frac{n-1}{2}(p-1)}  |\mathcal{V}(y+R,y)|^p \, \mathrm{d}y \label{lb estimate U nlin}
\end{align} for $z\geqslant R$, where we employed $ \mathcal{V}(y-R,y)=0$ that follows from the support condition \eqref{support condition aux fun}.
For $\mathcal{V}^{\mathrm{nlin}}$ the estimate from below on the characteristic line $t-z=R$ can be obtained in a similar way. For $z\geqslant R$ it holds
\begin{align}
\mathcal{V}^{\mathrm{nlin}}(R+z,z) 
&  \gtrsim \int_{R}^{z}  (R+y)^{-\frac{n-1}{2}(q-1)}  |\mathcal{U}(y+R,y)|^q \, \mathrm{d}y. \label{lb estimate V nlin}
\end{align}

Therefore, since $u_0,u_1,v_0,v_1$ and the kernel functions in the definitions of $\mathcal{U}^{\mathrm{lin}},\mathcal{V}^{\mathrm{lin}}$ are nonnegative,  for suitable positive constants $C,K$ depending on $n,p,q,R$ from \eqref{lb estimate U nlin} and \eqref{lb estimate V nlin} we have the iteration frame
\begin{align}
\mathcal{U}(R+z,z) & \geqslant C \int_{R}^{z} \mathrm{e}^{-\frac{b}{2}(z-y)} (R+y)^{-\frac{n-1}{2}(p-1)}  |\mathcal{V}(y+R,y)|^p \, \mathrm{d}y \, \quad \mbox{for} \ z\geqslant R, \label{IF 1} \\
\mathcal{V}(R+z,z) 
&  \geqslant K \int_{R}^{z}  (R+y)^{-\frac{n-1}{2}(q-1)}  |\mathcal{U}(y+R,y)|^q \, \mathrm{d}y \qquad  \qquad \quad \, \mbox{for} \ z\geqslant R. \label{IF 2} 
\end{align}

In order to start the iteration procedure, we need a first lower bound estimate for $\mathcal{V}(R+z,z)$. Since $v_0$ is nonnegative (and so is $\mathcal{V}_0$), from the definition of $\mathcal{V}^{\mathrm{lin}}$ we get immediately
\begin{align*}
\mathcal{V}^{\mathrm{lin}}(t,z)\geqslant  \frac{\varepsilon}{2} \int_{z-t}^{z+t}\mathcal{V}_1(y) \, \mathrm{d}y.
\end{align*} On the characteristic line $t-z=R$ for $z\geqslant R$ it results $[-R,R]\subset [z-t,z+t]$ , thus, from the support condition \eqref{support condition aux fun CD} we obtain
\begin{align}\label{lb V lin}
\mathcal{V}^{\mathrm{lin}}(R+z,z)\geqslant  \frac{\varepsilon}{2} \int_{\mathbb{R}} \mathcal{V}_1(y) \, \mathrm{d}y =   \frac{\varepsilon}{2} \int_{\mathbb{R}} \int_{\mathbb{R}^{n-1}} v_1(y,w) \, \mathrm{d}w \, \mathrm{d}y = \frac12 \| v_1\|_{L^1(\mathbb{R}^n)} \, \varepsilon, 
\end{align} where we used Fubini's theorem and the nonnegativity of $v_1$.

\begin{remark} Let us point out that for $\mathcal{U}^{\mathrm{lin}}$ we may derive only lower bounds that decay exponentially. Namely, since $\mathrm{I}_0(s)\geqslant 0$ and $\mathrm{I}_1(s)\geqslant \frac{s}{2}$ for $s\geqslant 0$ (the estimate from below for $\mathrm{I}_1$ is a straightforward consequence of the Maclaurin series expansion), and we assumed $u_0,u_1\geqslant 0$, from the definition of $\mathcal{U}^{\mathrm{lin}}$ for $z\geqslant R$ we have
\begin{align*}
\mathcal{U}^{\mathrm{lin}}(R+z,z)\geqslant \frac{\varepsilon}{2} \big\| u_1+\tfrac b2 u_0\big\|_{L^1(\mathbb{R}^n)} \, \mathrm{e}^{-\frac{b}{2}t} + \frac{\mu^2 \varepsilon}{4} \| u_0\|_{L^1(\mathbb{R}^n)} \, t \, \mathrm{e}^{-\frac{b}{2}t}.
\end{align*} Unfortunately, combining the previous exponential lower bound for $\mathcal{U}$ with the iteration frame \eqref{IF 1}-\eqref{IF 2} we are not able to get a sequence of lower bound estimates for $\mathcal{U}(R+z,z)$ whose lower bound diverges as $j\to\infty$ for $t$ above a certain $\varepsilon$-dependent threshold ($j$ denotes here the index in the sequence of lower bounds). In other words, an exponentially decaying lower bound for $\mathcal{U}$ does not allow us to derive a blow-up result for \eqref{Nakao der typ}.
\end{remark}

Since the nonlinear term in the second equation in \eqref{Nakao der typ} is nonnegative, from \eqref{lb V lin} it follows
\begin{align}\label{1st lb mathcal V}
\mathcal{V}^{}(R+z,z)\geqslant   M\varepsilon
\end{align} for $z\geqslant R$, where $M\doteq \frac12 \| v_1\|_{L^1(\mathbb{R}^n)}$.

We can start now the iteration argument to get a sequence of lower bound estimates for $\mathcal{V}(R+z,z)$. Since in \eqref{IF 1} it is present an exponential factor we need to use a slicing procedure when shrinking the domain of integration. The idea to shrink the domain of integration and cut intervals smaller and smaller on each step (i.e. the \emph{slicing procedure}) was introduced for the first time in \cite{AKT00}. Hence, in the series of papers \cite{ChenPal20,ChenPal20der} it was developed a slicing procedure associated with an increasing exponential function. Later, this method has been applied to study the blow-up dynamic of several semilinear weakly coupled systems (cf. \cite{ChenRei21,Chen21NA,Chen21MMAS}).

We shall consider separately the treatment of the subcritical case $\theta(n,p,q)>0$ from the treatment of the critical case $\theta(n,p,q)=0$.

\subsection{Subcritical case} \label{Subsection subcrit case}

In this section we focus on the subcritical case $\theta(n,p,q)>0$.
Let us introduce the parameters that individuate the slicing procedure, namely, the sequences of positive reals $\{\ell_j\}_{j\in\mathbb{N}}$, $\{L_j\}_{j\in\mathbb{N}}$ defined as follows:
\begin{align}
\ell_0 & \doteq \max\left\{\frac{2}{bR} ,1 \right\}\quad \mbox{and} \quad \ell_j \doteq 1+(pq)^{-j} \label{def ell j} \quad \mbox{for any} \  j\in\mathbb{N} \smallsetminus \{0\}, \\
L_j & \doteq \prod_{k=0}^j \ell_k \label{def L j} \quad \mbox{for any} \  j\in\mathbb{N}. 
\end{align} 
We emphasize that 
\begin{align}\label{def L}
L\doteq \lim_{j\to \infty} L_j= \prod_{k=0}^\infty \ell_k \in \mathbb{R}
\end{align} and, moreover, since $\ell_j>1$ for any $j\in\mathbb{N}\smallsetminus \{0\}$, it results $L_j\uparrow L$ as $j\to \infty$.

Our next goal is to prove
\begin{align}\label{lb mathcal V j} 
\mathcal{V}(R+z,z)\geqslant C_j (R+z)^{-\alpha_j}(z-L_j R)^{\beta_j} \qquad \mbox{for} \ z\geqslant L_j R \ \ \mbox{and any}  \  j\in\mathbb{N},
\end{align} where $\{C_j\}_{j\in\mathbb{N}}$, $\{\alpha_j\}_{j\in\mathbb{N}}$ and $\{\beta_j\}_{j\in\mathbb{N}}$ are sequences of nonnegative real numbers that we shall determine iteratively.
Clearly, due to \eqref{1st lb mathcal V}, \eqref{lb mathcal V j} for $j=0$ holds true by setting $C_0\doteq M\varepsilon$ and $\alpha_0\doteq 0,\beta_0\doteq 0$. Next we prove the inductive step. We assume that \eqref{lb mathcal V j} is satisfied for some $j\geqslant 0$ and we will prove it for $j+1$. Plugging \eqref{lb mathcal V j} in \eqref{IF 1}, for $z\geqslant L_jR$  we get
\begin{align*}
\mathcal{U}(R+z,z) & \geqslant C \int_{L_j R}^z \mathrm{e}^{-\frac{b}{2}(z-y)} (R+y)^{-\frac{n-1}{2}(p-1)}  |\mathcal{V}(y+R,y)|^p \, \mathrm{d}y \\
& \geqslant CC_j^p \int_{L_j R}^z \mathrm{e}^{-\frac{b}{2}(z-y)} (R+y)^{-\frac{n-1}{2}(p-1)-\alpha_j p}  (y-L_j R)^{\beta_j p} \, \mathrm{d}y \\
& \geqslant CC_j^p (R+z)^{-\frac{n-1}{2}(p-1)-\alpha_j p}  \int_{L_j R}^z \mathrm{e}^{-\frac{b}{2}(z-y)}  (y-L_j R)^{\beta_j p} \, \mathrm{d}y.
\end{align*} Thus, if we consider $z\geqslant L_{j+1}R$ then $[z/\ell_{j+1},z]\subset [L_j R,z]$. Therefore, shrinking the domain of integration in the previous inequality, for $z\geqslant L_{j+1}R$ we have
\begin{align}
\mathcal{U}(R+z,z) &  \geqslant CC_j^p (R+z)^{-\frac{n-1}{2}(p-1)-\alpha_j p}  \int_{z/\ell_{j+1}}^z \mathrm{e}^{-\frac{b}{2}(z-y)}  (y-L_j R)^{\beta_j p} \, \mathrm{d}y \notag \\
& \geqslant CC_j^p (R+z)^{-\frac{n-1}{2}(p-1)-\alpha_j p} \left(\tfrac{z}{\ell_{j+1}}-L_j R\right)^{\beta_j p} \int_{z/\ell_{j+1}}^z \mathrm{e}^{-\frac{b}{2}(z-y)}   \, \mathrm{d}y \notag \\
& = 2b^{-1}CC_j^p \ell_{j+1}^{-\beta_j p} (R+z)^{-\frac{n-1}{2}(p-1)-\alpha_j p} (z-L_j \ell_{j+1} R)^{\beta_j p} \left(1-\mathrm{e}^{-\frac{b}{2}(1-1/\ell_{j+1})z}\right). \label{lb mathcal U intermediate}
\end{align} Let us estimate from below the factor on the right hand-side of the previous chain of inequalities that contains the exponential term. Then, for $z\geqslant L_{j+1}R$ it holds
\begin{align}
1-\mathrm{e}^{-\frac{b}{2}(1-1/\ell_{j+1})z} & \geqslant 1-\mathrm{e}^{-\frac{b}{2}RL_{j+1}(1-1/\ell_{j+1})} = 1-\mathrm{e}^{-\frac{b}{2}RL_{j}(\ell_{j+1}-1)}\geqslant 1-\mathrm{e}^{-\frac{b}{2}R\ell_{0}(\ell_{j+1}-1)} \notag \\
&  \geqslant  1-\mathrm{e}^{-(\ell_{j+1}-1)} \geqslant 1-\left(1-(\ell_{j+1}-1)+\tfrac 12 (\ell_{j+1}-1)^2\right) = (\ell_{j+1}-1)\left(1-\tfrac 12 (\ell_{j+1}-1)\right) \notag \\
&= (pq)^{-2(j+1)}\left((pq)^{j+1}-\tfrac 12\right)\geqslant (pq)^{-2(j+1)}\left((pq)-\tfrac 12\right). \label{lb est exp term}
\end{align} Combining \eqref{lb mathcal U intermediate} and \eqref{lb est exp term}, for $z\geqslant L_{j+1} R$ we arrive at
\begin{align*}
\mathcal{U}(R+z,z) &  \geqslant  (2pq-1)b^{-1}CC_j^p \ell_{j+1}^{-\beta_j p} (pq)^{-2(j+1)} (R+z)^{-\frac{n-1}{2}(p-1)-\alpha_j p} (z-L_{j+1} R)^{\beta_j p}.
\end{align*} Plugging the previous upper bound for $\mathcal{U}(R+z,z)$ in \eqref{IF 2}, for $z\geqslant L_{j+1} R$ we get
\begin{align*}
\mathcal{V}(R+z,z) 
&  \geqslant K \int_{ L_{j+1} R}^{z}  (R+y)^{-\frac{n-1}{2}(q-1)}  |\mathcal{U}(y+R,y)|^q \, \mathrm{d}y \\
&  \geqslant \frac{KC^q (2pq-1)^q b^{-q} C_j^{pq} }{\ell_{j+1}^{\beta_j pq} (pq)^{2q(j+1)} }\int_{ L_{j+1} R}^{z}  (R+y)^{-\frac{n-1}{2}(pq-1)-\alpha_j pq} (y-L_{j+1} R)^{\beta_j pq} \, \mathrm{d}y \\
&  \geqslant \frac{KC^q (2pq-1)^q b^{-q} C_j^{pq} }{\ell_{j+1}^{\beta_j pq} (pq)^{2q(j+1)} }  (R+z)^{-\frac{n-1}{2}(pq-1)-\alpha_j pq}\int_{ L_{j+1} R}^{z}  (y-L_{j+1} R)^{\beta_j pq} \, \mathrm{d}y \\
&  = \frac{KC^q (2pq-1)^q b^{-q} C_j^{pq} }{\ell_{j+1}^{\beta_j pq} (pq)^{2q(j+1)} (\beta_j pq+1)}  (R+z)^{-\frac{n-1}{2}(pq-1)-\alpha_j pq}  (z-L_{j+1} R)^{\beta_j pq+1}. 
\end{align*} Thus, we proved \eqref{lb mathcal V j} for $j+1$ with
\begin{align}
C_{j+1} &\doteq \frac{KC^q (2pq-1)^q b^{-q} C_j^{pq} }{\ell_{j+1}^{\beta_j pq} (pq)^{2q(j+1)} (\beta_j pq+1)}, \label{def C j+1} \\
\alpha_{j+1} & \doteq \frac{n-1}{2}(pq-1)+pq \alpha_j, \qquad \beta_{j+1}\doteq 1+pq \beta_j. \label{def alpha beta j+1}
\end{align}
The next step is to determine a suitable lower bound for $C_j$, that will be easier to handle. First we derive an explicit representation for $\alpha_j$ and $\beta_j$. By using recursively \eqref{def alpha beta j+1}, we have
\begin{align}
\alpha_j & = \tfrac{n-1}{2}(pq-1)+pq \alpha_{j-1}= \cdots= (pq)^j\alpha_0+\tfrac{n-1}{2}(pq-1)\sum_{k=0}^{j-1}(pq)^k  =\tfrac{n-1}{2}((pq)^j-1), \label{representation alpha j}\\
\beta_j &= 1+pq\beta_{j-1}= \cdots = (pq)^j\beta_0+\sum_{k=0}^{j-1}(pq)^k = \tfrac{(pq)^j-1}{pq-1}.  \label{representation beta j}
\end{align} Therefore,
\begin{align*}
(\beta_{j-1}pq+1)^{-1}=\beta_j^{-1}\geqslant (pq-1) (pq)^{-j}.
\end{align*} Moreover, since 
\begin{align*}
\lim_{j\to\infty} \ell^{\beta_{j-1}pq}_j=\lim_{j\to\infty} \exp\left(\frac{(pq)^j-pq}{pq-1}\ln\big(1+(pq)^{-j}\big)\right)=\mathrm{e}^{1/(pq-1)}, 
\end{align*} there exists $N=N(p,q,b,R)>0$ such that $\ell_j^{-\beta_{j-1}pq}>N$ for any $j\in\mathbb{N}$. Consequently,
\begin{align}\label{lb Cj}
C_j=\frac{KC^q (2pq-1)^q b^{-q} C_{j-1}^{pq} }{\ell_{j}^{\beta_{j-1} pq} (pq)^{2qj} (\beta_{j-1} pq+1)}\geqslant D (pq)^{-(2q+1)j}C_{j-1}^{pq}
\end{align} for any $j\in\mathbb{N}$, where $D\doteq KC^q N (2pq-1)^q(pq-1) b^{-q}$. Applying the logarithmic function to both sides of \eqref{lb Cj} and using iteratively the resulting inequality, we have
\begin{align*}
\ln C_j &\geqslant pq\ln C_{j-1} -(2q+1)j\ln (pq)+\ln D \\
& \geqslant  (pq)^2\ln C_{j-2} -(2q+1)\ln (pq)(j+(j-1)pq) +(1+pq)\ln D \\
& \geqslant \cdots\geqslant (pq)^j\ln C_{0} -(2q+1)\ln (pq)\sum_{k=0}^{j-1} (j-k)(pq)^k +\ln D \sum_{k=0}^{j-1} (pq)^k.
\end{align*} Using the identities
\begin{align}\label{summation identities}
\sum_{k=0}^{j-1} (j-k)(pq)^k =\frac{1}{pq-1}\left(\frac{(pq)^{j+1}-pq}{pq-1}-j\right), \qquad \sum_{k=0}^{j-1} (pq)^k= \frac{(pq)^j-1}{pq-1},
\end{align} it results
\begin{align*}
\ln C_j & \geqslant (pq)^j\left(\ln C_{0} -\frac{(2q+1)pq\ln(pq)}{(pq-1)^2}+\frac{\ln D}{pq-1}\right)+\frac{(2q+1)pq\ln(pq)}{(pq-1)^2}+\frac{(2q+1)\ln(pq)}{pq-1}j-\frac{\ln D}{pq-1}.
\end{align*} Let us denote by $j_0=j_0(n,b,p,q,R)$ the smallest nonnegative integer such that
\begin{align*}
j_0\geqslant \frac{\ln D}{(2q+1)\ln(pq)}-\frac{pq}{pq-1}.
\end{align*} Then, for any $j\geqslant j_0$
\begin{align} \label{lb Cj final}
\ln C_j & \geqslant (pq)^j\left(\ln C_{0} -\frac{(2q+1)pq\ln(pq)}{(pq-1)^2}+\frac{\ln D}{pq-1}\right)=(pq)^j\ln(E\varepsilon),
\end{align} where $E\doteq M(pq)^{-(2q+1)(pq)/(pq-1)^2}D^{1/(pq-1)}$.
Combining \eqref{lb mathcal V j}, \eqref{representation alpha j}, \eqref{representation beta j} and \eqref{lb Cj final}, for $j\geqslant j_0$ and $z\geqslant LR$ we find
\begin{align*}
\mathcal{V}(R+z,z) &\geqslant \exp\left((pq)^j\ln(E\varepsilon)\right) (R+z)^{-\frac{n-1}{2}((pq)^j-1)}(z-L R)^{\frac{(pq)^j-1}{pq-1}} \\
& = \exp\left((pq)^j\left(\ln(E\varepsilon)-\tfrac{n-1}{2}\ln(R+z)+\tfrac{1}{pq-1}\ln(z-LR)\right)\right) (R+z)^{\frac{n-1}{2}}(z-L R)^{-\frac{1}{pq-1}}.
\end{align*} Equivalently, for $t\geqslant (L+1)R$ and $j\geqslant j_0$ it holds
\begin{align*}
\mathcal{V}(t,t-R) &\geqslant \exp\left((pq)^j\left(\ln(E\varepsilon)-\tfrac{n-1}{2}\ln t +\tfrac{1}{pq-1}\ln(t-(L+1)R)\right)\right) t^{\frac{n-1}{2}}(t-(L+1) R)^{-\frac{1}{pq-1}}.
\end{align*} For $t\geqslant 2(L+1)R$ we can estimate $\ln (t-(L+1)R)\geqslant \ln t -\ln 2$. Consequently, for $t\geqslant 2(L+1)R$ and $j\geqslant j_0$ we have
\begin{align}
\mathcal{V}(t,t-R) &\geqslant \exp\left((pq)^j\left(\ln(E\varepsilon)+(\tfrac{1}{pq-1}-\tfrac{n-1}{2})\ln t -\tfrac{1}{pq-1} \ln 2 \right)\right) t^{\frac{n-1}{2}}(t-(L+1) R)^{-\frac{1}{pq-1}} \notag \\
& =\exp\left((pq)^j\left(\ln\left(E_1\varepsilon t^{\theta(n,p,q)}\right)\right)\right) t^{\frac{n-1}{2}}(t-(L+1) R)^{-\frac{1}{pq-1}}, \label{lb mathcal V final}
\end{align} where $E_1\doteq 2^{-1/(pq-1)}E$ and $\theta$ is defined in \eqref{def theta}.

Let us fix $\varepsilon_0=\varepsilon_0(n,p,q,b,R,v_1)>0$ sufficiently small so that $$\varepsilon_0\leqslant E_1^{-1} (2(L+1)R)^{-\theta(n,p,q)}.$$ Then, for any $\varepsilon\in(0,\varepsilon_0]$ and any $t\geqslant (E_1 \varepsilon)^{-\theta(n,p,q)^{-1}}$ we have
\begin{align*}
t\geqslant 2(L+1)R \quad \mbox{and} \quad \ln\left(E_1\varepsilon t^{\theta(n,p,q)}\right)>0,
\end{align*} so letting $j\to\infty$ in \eqref{lb mathcal V final} we see that $\mathcal{V}(t,t-R)$ is not finite. 

Summarizing we proved that $(u,v)$ blows up in finite time and we established the upper bound estimate in \eqref{upper bound lifespan thm}.

\subsection{Critical case}
In this section we study the blow-up result  in the critical case $\theta(n,p,q)=0$. In this case it is more convenient to rewrite the iteration frame as follows 
\begin{align}
\mathcal{U}(R+z,z) & \geqslant C \int_{R}^{z} \mathrm{e}^{-\frac{b}{2}(z-y)} y^{-\frac{n-1}{2}(p-1)}  |\mathcal{V}(R+y,y)|^p \, \mathrm{d}y \, \quad \mbox{for} \ z\geqslant R, \label{IF 1 crit} \\
\mathcal{V}(R+z,z) 
&  \geqslant K \int_{R}^{z}  y^{-\frac{n-1}{2}(q-1)}  |\mathcal{U}(R+y,y)|^q \, \mathrm{d}y \qquad  \qquad \quad \, \mbox{for} \ z\geqslant R. \label{IF 2 crit} 
\end{align} Note that, for the sake of simplicity, we kept the same notations for the multiplicative constants as in Section \ref{Subsection subcrit case}.

The main difference in comparison to the subcritical case consists in the choice of the parameters characterizing the slicing procedure. We introduce the sequence $\{\Lambda_j\}_{j\in\mathbb{N}}$, where $$\Lambda_j\doteq 1+\frac{4}{bR}\big(2-2^{-j}\big) \quad \mbox{for any}  \ j\in\mathbb{N}.$$ The sequence $\{\Lambda_j\}_{j\in\mathbb{N}}$ is strictly increasing and bounded and, clearly, $\Lambda_j \uparrow \Lambda \doteq 1+8/(bR)$ as $j\to \infty$. We shall employ this sequence when applying the slicing procedure. 

The next step is to prove the sequence of lower bound estimates
\begin{align}\label{lb mathcal V j crit} 
\mathcal{V}(R+z,z)\geqslant K_j \left(\ln\left(\frac{z}{\Lambda_j R}\right)\right)^{\gamma_j} \qquad \mbox{for} \ z\geqslant \Lambda_j R \ \ \mbox{and any}  \  j\in\mathbb{N},
\end{align} where $\{K_j\}_{j\in\mathbb{N}}$ and $\{\gamma_j\}_{j\in\mathbb{N}}$ are sequences of nonnegative real numbers to be determined iteratively throughout the proof. 

We remark that \eqref{1st lb mathcal V} implies the validity of \eqref{lb mathcal V j crit} for $j=0$, provided that $K_0=M\varepsilon$ and $\gamma_0=0$. In order to establish \eqref{lb mathcal V j crit} for any $j\in\mathbb{N}$ it remains to demonstrate the inductive step. Let us assume that \eqref{lb mathcal V j crit} is fulfilled for some $j\in\mathbb{N}$, then, we have to prove that \eqref{lb mathcal V j crit} is satisfied also for $j+1$. Plugging \eqref{lb mathcal V j crit} in \eqref{IF 1 crit}, for $z\geqslant \Lambda_{j+1} R$ we find
\begin{align*}
\mathcal{U}(R+z,z) & \geqslant C \int_{\Lambda_j R}^{z} \mathrm{e}^{-\frac{b}{2}(z-y)} y^{-\frac{n-1}{2}(p-1)}  |\mathcal{V}(R+y,y)|^p \, \mathrm{d}y \\
& \geqslant C  K_j^p \int_{\Lambda_j R}^{z} \mathrm{e}^{-\frac{b}{2}(z-y)} y^{-\frac{n-1}{2}(p-1)}  \left(\ln\left(\frac{y}{\Lambda_j R}\right)\right)^{\gamma_j p} \, \mathrm{d}y  \\
& \geqslant C  K_j^p z^{-\frac{n-1}{2}(p-1)} \int_{\tfrac{\Lambda_j z}{\Lambda_{j+1}} }^{z} \mathrm{e}^{-\frac{b}{2}(z-y)}   \left(\ln\left(\frac{y}{\Lambda_j R}\right)\right)^{\gamma_j p} \, \mathrm{d}y \\
& \geqslant C  K_j^p z^{-\frac{n-1}{2}(p-1)} \left(\ln\left(\frac{z}{\Lambda_{j+1} R}\right)\right)^{\gamma_j p} \int_{\tfrac{\Lambda_j z}{\Lambda_{j+1}} }^{z} \mathrm{e}^{-\frac{b}{2}(z-y)}    \, \mathrm{d}y  \\
& =2 b^{-1}C  K_j^p z^{-\frac{n-1}{2}(p-1)} \left(\ln\left(\frac{z}{\Lambda_{j+1} R}\right)\right)^{\gamma_j p} \left(1- \mathrm{e}^{-\frac{b}{2}(\Lambda_{j+1}-\Lambda_j)\tfrac{z}{\Lambda_{j+1}}}\right),
\end{align*} where in the third step we used $[\Lambda_j  R,z]\supset \big[\frac{\Lambda_j z}{\Lambda_{j+1}},z\big]$. Since  $\{\Lambda_j\}_{j\in\mathbb{N}}$ is an increasing sequence, for  $z\geqslant \Lambda_{j+1} R$ we may estimate 
\begin{align*}
1- \mathrm{e}^{-\frac{b}{2}(\Lambda_{j+1}-\Lambda_j)\tfrac{z}{\Lambda_{j+1}}} & \geqslant 1- \mathrm{e}^{-\frac{bR}{2}(\Lambda_{j+1}-\Lambda_j)} \\
& \geqslant \frac{bR}{2}(\Lambda_{j+1}-\Lambda_j)\left(1-\frac{bR}{4}(\Lambda_{j+1}-\Lambda_j)\right) \\
&= 2^{-(2j+1)}(2^{j+1}-1)\geqslant 2^{-(2j+1)},
\end{align*}  where in the second inequality we used the elementary inequality $\mathrm{e}^{-s}\leqslant 1-s+\frac{s^2}{2}$ for any $s\geqslant 0$. Thus, for $z\geqslant \Lambda_{j+1} R$ we showed that 
\begin{align*}
\mathcal{U}(R+z,z) & \geqslant b^{-1}C  \, 2^{-2j}  K_j^p z^{-\frac{n-1}{2}(p-1)} \left(\ln\left(\frac{z}{\Lambda_{j+1} R}\right)\right)^{\gamma_j p} .
\end{align*} Using the last lower bound estimate for $\mathcal{U}(R+z,z)$ in \eqref{IF 2 crit} and the critical condition $\theta(n,p,q)=0$,  for $z\geqslant \Lambda_{j+1} R$ we get
\begin{align*}
\mathcal{V}(R+z,z) 
&  \geqslant K \int_{\Lambda_{j+1} R}^{z}  y^{-\frac{n-1}{2}(q-1)}  |\mathcal{U}(R+y,y)|^q \, \mathrm{d}y \\
&  \geqslant K b^{-q} C^q  \, 2^{-2qj}  K_j^{pq}  \int_{\Lambda_{j+1} R}^{z}  y^{-\frac{n-1}{2}(pq-1)} \left(\ln\left(\frac{y}{\Lambda_{j+1} R}\right)\right)^{\gamma_j pq} \, \mathrm{d}y \\
&  = K b^{-q} C^q  \, 2^{-2qj}  K_j^{pq}  \int_{\Lambda_{j+1} R}^{z}  y^{-1} \left(\ln\left(\frac{y}{\Lambda_{j+1} R}\right)\right)^{\gamma_j pq} \, \mathrm{d}y \\
&  = K b^{-q} C^q  \, 2^{-2qj}  K_j^{pq} (\gamma_j pq+1)^{-1} \left(\ln\left(\frac{z}{\Lambda_{j+1} R}\right)\right)^{\gamma_j pq+1},
\end{align*} which is exactly \eqref{lb mathcal V j crit} for $j+1$, provided that we set 
\begin{align}
K_{j+1} &\doteq K b^{-q} C^q  \, 2^{-2qj}  K_j^{pq} (\gamma_j pq+1)^{-1}, \label{def K j+1} \\
\gamma_{j+1} &\doteq \gamma_j pq+1. \label{def gamma j+1}
\end{align}
By applying iteratively \eqref{def gamma j+1} and $\gamma_0=0$, we obtain
\begin{align} \label{def gamma j}
\gamma_j=  1+pq\gamma_{j-1}= \cdots = (pq)^j\gamma_0+\sum_{k=0}^{j-1}(pq)^k = \tfrac{(pq)^j-1}{pq-1}.  
\end{align} Next we determine a lower bound estimate for the constant $K_j$. From the previous representation for $\gamma_j$ it follows
\begin{align*}
K_{j} & = K b^{-q} C^q  \, 2^{-2q(j-1)}  K_{j-1}^{pq} (\gamma_{j-1} pq+1)^{-1} = K b^{-q} C^q  \, 2^{-2q(j-1)}  K_{j-1}^{pq} \gamma_{j}^{-1} \\ &\geqslant 2^{2q} K b^{-q} C^q (pq-1) \, 2^{-2qj} (pq)^{-j} K_{j-1}^{pq} = \tilde{D} (2^{2q}pq)^{-j}  K_{j-1}^{pq}, 
\end{align*} where $\tilde{D}\doteq 2^{2q} K b^{-q} C^q (pq-1)$. Applying the logarithmic function to both sides of the inequality $K_{j}\geqslant \tilde{D} (2^{2q}pq)^{-j}  K_{j-1}^{pq}$ and using in an iterative way the obtained inequality, we obtain
\begin{align*}
\ln K_j &\geqslant pq\ln K_{j-1} -j\ln (2^{2q}pq)+\ln \tilde{D} \\
& \geqslant  (pq)^2\ln K_{j-2} -\ln (2^{2q}pq)(j+(j-1)pq) +(1+pq)\ln \tilde{D} \\
& \geqslant \cdots\geqslant (pq)^j\ln K_{0} -\ln (2^{2q}pq)\sum_{k=0}^{j-1} (j-k)(pq)^k +\ln \tilde{D} \sum_{k=0}^{j-1} (pq)^k\\
& =(pq)^j \left(\ln (M\varepsilon) -\frac{(pq)\ln (2^qpq) }{(pq-1)^2}+\frac{\ln \tilde{D}}{pq-1}\right)+\frac{(pq)\ln (2^q pq)}{(pq-1)^2}+\frac{\ln (2^q pq)}{pq-1}j-\frac{\ln \tilde{D}}{pq-1},
\end{align*} where in the last step we applied the identities in \eqref{summation identities}. If we denote by $j_1=j_1(n,b,p,q,R)$ the smallest nonnegative integer number such that $$j_1\geqslant \frac{\ln \tilde{D}}{\ln (2^q pq)}-\frac{pq}{pq-1},$$ then, for any $j\geqslant j_1$ it results
\begin{align}\label{lb Kj}
\ln K_j &\geqslant (pq)^j \left(\ln (M\varepsilon) -\frac{(pq)\ln (2^qpq) }{(pq-1)^2}+\frac{\ln \tilde{D}}{pq-1}\right) = (pq)^j \ln (\tilde{E}\varepsilon),
\end{align} where $\tilde{E}\doteq  M (2^q pq)^{-(pq)/(pq-1)^2}\tilde{D}^{1/(pq-1)}$. 

Combining \eqref{lb mathcal V j crit}, \eqref{def gamma j} and \eqref{lb Kj}, for $j\geqslant j_1$ and $z\geqslant \Lambda R$ we find
\begin{align*}
\mathcal{V}(R+z,z) & \geqslant \exp\left((pq)^j\ln (\tilde{E}\varepsilon)\right) \left(\ln\left(\frac{z}{\Lambda R}\right)\right)^{\tfrac{(pq)^j-1}{pq-1}}  \\
&=  \exp\left((pq)^j\left(\ln \left(\tilde{E}\varepsilon\ln\left(\frac{z}{\Lambda R}\right)^{\frac{1}{pq-1}}\right)\right)\right) \left(\ln\left(\frac{z}{\Lambda R}\right)\right)^{-\tfrac{1}{pq-1}}  
\end{align*} Therefore, for $t\geqslant (\Lambda+1)R$ 
and for any $j\geqslant j_1$ we have
\begin{align}
\mathcal{V}(t,t-R) & \geqslant   \exp\left((pq)^j\left(\ln \left(\tilde{E}\varepsilon\ln\left(\frac{t-R}{\Lambda R}\right)^{\frac{1}{pq-1}}\right)\right)\right) \left(\ln\left(\frac{t-R}{\Lambda R}\right)\right)^{-\tfrac{1}{pq-1}} \notag \\
& \geqslant   \exp\left((pq)^j\left(\ln \left(\tilde{E}\varepsilon\ln\left(\frac{t}{2\Lambda R}\right)^{\frac{1}{pq-1}}\right)\right)\right) \left(\ln\left(\frac{t-R}{\Lambda R}\right)\right)^{-\tfrac{1}{pq-1}} . \label{last lb mathcal V crit}
\end{align}
Let us fix $\varepsilon_0=\varepsilon_0(n,b,p,q,R)>0$ sufficiently small so that $\varepsilon_0^{-(pq-1)}\geqslant \tilde{E}^{pq-1}\ln \frac{\Lambda+1}{2\Lambda}$. Then, for any $\varepsilon\in (0,\varepsilon_0]$ and any $t> (2\Lambda R)\exp \big((\tilde{E}\varepsilon)^{-(pq-1)}\big)$ the following inequalities are satisfied
\begin{align*}
t\geqslant (\Lambda+1)R \quad  \mbox{and} \quad \ln \left(\tilde{E}\varepsilon\ln\left(\frac{t}{2\Lambda R}\right)^{\frac{1}{pq-1}}\right)>1,
\end{align*} thus, letting $j\to \infty$ in \eqref{last lb mathcal V crit} the lower bound for $\mathcal{V}(t,t-R)$ diverges. Consequently, $\mathcal{V}(t,t-R)$ cannot be finite. So, we proved that $v$ blows up in finite time and we found as byproduct of the iteration procedure the upper bound estimate $$T(\varepsilon)\leqslant \exp\left(\tilde{E}_1 \varepsilon^{-(pq-1)}\right)$$ for the lifespan of the solution $(u,v)$ for any $\varepsilon\in(0,\varepsilon_0]$, where $\tilde{E}_1>0$ is a suitable constant depending on $n,b,p,q,R,v_1$.


\section*{Acknowledgments}

A. Palmieri is supported by the\emph{ Japan Society for the Promotion of Science} (JSPS) – JSPS Postdoctoral Fellowship for Research in Japan (Short-term) (PE20003) – and is member of the \emph{Gruppo Nazionale per L’Analisi Matematica, la Probabilità e le loro Applicazioni} (GNAMPA) of the \emph{Instituto Nazionale di Alta Matematica} (INdAM). H.Takamura is partially supported by the Grant-in-Aid for Scientific Research (B) (No.18H01132), \emph{Japan Society for the Promotion of Science}.

\end{document}